\documentclass{article}

\usepackage{fullpage}
\usepackage{verbatim}
\usepackage{amsmath}
\usepackage{amssymb}
\usepackage{amsfonts, pstricks}
\usepackage{epsf}
\usepackage[dvips]{graphicx}
\usepackage[comma, square]{natbib}
\usepackage[latin1]{inputenc}

\newtheorem{guil}{Observation}

\newtheorem{schr}[guil]{Observation}
\newtheorem{guil2}[guil]{Theorem}
\newtheorem{guil3}[guil]{Corollary}
\newtheorem{alt0}[guil]{Theorem}

\newtheorem{wind}[guil]{Theorem}

\newtheorem{permd}[guil]{Theorem}
\newtheorem{sepenum}[guil]{Observation}
\newtheorem{forpat}[guil]{Observation}
\newtheorem{permq}[guil]{Corollary}
\newtheorem{sep3char}[guil]{Theorem}

\begin{document}

\renewcommand{\thesubsection}{\arabic{subsection}}

\author{Andrei Asinowski
\footnote{Caesarea Rothschild Institute, University of Haifa, Haifa 31905, Israel. E-mail \texttt{andrei@cri.haifa.ac.il}.}
\and Toufik Mansour
\footnote{Department of Mathematics, University of Haifa, Haifa 31905, Israel. E-mail \texttt{toufik@math.haifa.ac.il}. }}

\title{Separable $d$-permutations and guillotine partitions}

\date{}

\maketitle

\begin{abstract}
We characterize separable multidimensional permutations in terms of forbidden patterns and enumerate them by means of generating function, recursive formula and explicit formula.
We find a connection between multidimensional permutations and guillotine partitions of a box.
In particular, a bijection between $d$-dimensional permutations and guillotine partitions of a $2^{d-1}$-dimensional box is constructed.
We also study enumerating problems related to guillotine partitions under certain restrictions revealing
connections to other combinatorial structures.
This allows us to obtain results on avoided patterns in permutations.
\end{abstract}

\textsc{AMS 2000 Subject Classification:} Primary 05A05, 05A15; Secondary 05C30, 68R05.

\textsc{Keywords:} $d$-permutations, separable permutations, patterns in permutations, guillotine partitions, binary trees, Schr\"oder paths.

\newcommand{\alter}{\begin{pspicture}(0,.1)(.5,.5)
\psline(0,0)(.1333,0)(.1333,.4)(.1333,0)(.2666,0)(.2666,.4)(.2666,0)(.4,0)(.4,.4)(0,.4)(0,0)
\end{pspicture}}

\newcommand{\winaa}{\begin{pspicture}(0,.1)(.5,.5)
\psline(0,0)(.2,0)(.2,.4)(.2,0)(.4,0)(.4,.16)(.2,.16)
(.4,.16)(.4,.4)(0,.4)(0,.24)(.2,.24)(0,.24)(0,0)
\end{pspicture}}

\newcommand{\winbb}{\begin{pspicture}(0,.1)(.5,.5)
\psline(0,0)(.2,0)(.2,.4)(.2,0)(.4,0)(.4,.24)(.2,
.24)(.4,.24)(.4,.4)(0,.4)(0,.16)(.2,.16)(0,.16)(0,0)
\end{pspicture}}

\newcommand{\wincc}{\begin{pspicture}(0,.1)(.5,.5)
\psline(0,0)(.133,0)(.133,.4)(.133,0)(.266,0)(.266,.4)(.266,0)(.4,0)(.4,.24)(.266,
.24)(.4,.24)(.4,.4)(0,.4)(0,.16)(.133,.16)(0,.16)(0,0)
\end{pspicture}}

\newcommand{\y}{\mathrm}

\section{Introduction}
In the first part of this paper we study the multidimensional generalization of separable permutations. Separable permutations form a well known class of permutations, they may be defined recursively as follows: a separable permutation is either a permutation of one element or a concatenation of two smaller separable permutations, upon an appropriate relabeling. A $d$-permutation is a sequence of $d$ permutations, the first of them being the natural order permutation $12\dots n$. The notion of separable permutations generalizes to that of separable $d$-permutations in a natural way. After formal definitions (Section~\ref{sec:perm_intro}), we find the generating function, a recursive formula and an explicit formula for the number of separable $d$-permutations of $\{1, 2, \dots, n\}$ (Section~\ref{sec:perm_enum}) and characterize them in terms of forbidden patterns (Section~\ref{sec:perm_forbid}).

The second part of the paper is devoted to guillotine partitions of a $d$-dimensional box, \emph{i.e.,}
recursive partitions of a $d$-dimensional box $B$ by axis-aligned hyperplanes. Guillotine partitions were introduced in $1980$ies, and they have numerous applications in computational geometry, computer graphics, \emph{etc}. Recently, Ackerman, Barequet, Pinter and Romik studied the enumerative issues related to guillotine partitions~\cite{ack3, ack2, ack}.
We observe that the generating function for the number of separable $d$-permutations is identical to the
generating function for the number of (structurally different) guillotine partitions of a $2^{d-1}$-dimensional box. Ackerman \emph{et al.} constructed a bijection between these sets in the case $d=2$.
In Section~\ref{sec:par_per} we generalize a version of their bijection to any $d$, and find a subclass of separable $d$-permutations which correspond to guillotine partitions of $q$-dimensional box where $q$ is not necessarily a power of $2$.

In Section~\ref{sec:restr} we deal with guillotine partitions with certain restrictions. In Sections~\ref{sec:bound} -- \ref{sec:win} we enumerate some classes of restricted guillotine partitions, in Section~\ref{sec:restr_par_per} we use these results and the correspondence between separable permutations and guillotine partitions for enumerating of permutations avoiding certain patterns.

\section{Separable $d$-permutations}
\label{sec:perm}

\subsection{Notation and convention}
\label{sec:perm_intro}

A \emph{$d$-permutation of $[n]=\{1, 2, \dots, n\}$} is a
sequence $P=(p_1, p_2, \dots, p_d)$ where each $p_i$ is a permutation of $[n]$ and $p_1$ is the natural-order permutation $12\dots n$. It may be represented as a $d\times n$ matrix (also denoted by $P$)
each row of which is a permutation of $[n]$, the first row being $12\dots n$. Thus, for $1 \leq i \leq d$, $p_i$ is a row of this matrix, and we shall denote by $P^{(j)}$ its $j$th column ($1 \leq j \leq n$). $P_{ij}$ will denote the $(i,j)$-entry of the matrix $P$.

A $d$-permutation may be represented geometrically as a point set in $\mathbb{R}^d$ -- a subset of size $n$ of the discrete cube $[n]^d$ such that each hyperplane $x_i = j$, $1 \leq i \leq d$, $1 \leq j \leq n$ contains precisely one point. We shall refer to this geometric representation as to the \emph{graph of $P$}. $P^{(j)}$ is the coordinate vector of the point whose first coordinate is $j$.

It is clear that there are $(n!)^{d-1}$ $d$-permutations of $[n]$.
\emph{Remark: We use one of two existing approaches to the notion of $d$-permutation. According to another approach, what we call $d$-permutation is called $(d-1)$-permutation. For example: According to our approach, usual permutations are $2$-permutations, according to the other they are $1$-permutations.}

A $d$-permutation $P$ of $[n]$ is \emph{separable} if
either $n=1$, or $n>1$ and there is a number $\ell$, $1 \leq \ell < n$, such that for each $1 \leq i \leq d$
we have
either $P_{ij_1}<P_{ij_2}$ or $P_{ij_1}>P_{ij_2}$
for all $1 \leq j_1 \leq \ell$, $\ell+1 \leq j_2 \leq n$,
and the two $d$-permutations obtained from $P$ by taking the first $\ell$ columns and
the last $n-\ell$ columns and order preserving relabeling them so that they become $d$-permutations of $[l]$ and $[n-\ell]$, respectively (these $d$-permutations will be called \emph{left} and \emph{right blocks of $P$ with respect to $\ell$} and denoted by $P_L^{\ell}$ and $P_R^{\ell}$), are themselves separable.
In this case we say that $P$ is \emph{separated between $\ell$ and $\ell+1$}.

A \emph{primary block structure} of a separable $d$-permutation $P$ (with $n\geq 2$), separated between $\ell$ and $\ell+1$, is a $d$-permutation $S=S(P)$ of $\{1, 2\}$ defined as follows:
\begin{quote}
if for each $1 \leq j_1 \leq \ell$, $\ell+1 \leq j_2 \leq n$ we have $P_{ij_1}<P_{ij_2}$, set $S_{i1}=1$, $S_{i2}=2$,\\
if for each $1 \leq j_1 \leq \ell$, $\ell+1 \leq j_2 \leq n$ we have $P_{ij_1}>P_{ij_2}$, set $S_{i1}=2$, $S_{i2}=1$.
\end{quote}

Note that in general, $\ell$ is not unique: a separable permutation may be separated in several places.
However, it is easy to see that the primary block structure is determined uniquely.

\vspace{2mm}

\noindent\textit{Example.} Consider

\[P=\left(
    \begin{array}{ccccc}
      1 & 2 & 3 & 4 & 5 \\ 
      2 & 1 & 3 & 4 & 5 \\ 
      5 & 4 & 3 & 1 & 2 \\
    \end{array}
  \right)
\]
It can be separated between $2$ and $3$ ($\ell = 2$), or between $3$ and $4$ ($\ell = 3$).
For $\ell = 2$,
\[P^2_L=\left(
    \begin{array}{cc}
      1 & 2 \\ 
      2 & 1 \\ 
      1 & 2 \\
    \end{array}
  \right), \ \
  P^2_R=\left(
    \begin{array}{ccc}
      1 & 2 & 3\\ 
      1 & 2 & 3\\ 
      3 & 1 & 2\\
    \end{array}
  \right);
\]
and for $\ell =3$,
\[P^3_L=\left(
    \begin{array}{ccc}
      1 & 2 & 3\\ 
      2 & 1 & 3\\ 
      3 & 2 & 1\\
    \end{array}
  \right), \ \
P^3_R=\left(
    \begin{array}{cc}
      1 & 2 \\ 
      1 & 2 \\ 
      1 & 2 \\
    \end{array}
  \right).
\]
The primary block structure of $P$ is
\[
S(P)=\left(
    \begin{array}{cc}
      1 & 2 \\ 
      1 & 2 \\ 
      2 & 1 \\
    \end{array}
  \right).
\]

Geometrically, the graph of a separated $d$-permutation of $[n]$, $n>1$, is obtained by placing graphs of two smaller separable $d$-permutations in opposite orthants, and appropriate relabeling.

For $d=2$, we get separable permutations.
It is well known that the number of separable permutations of $[n]$ is the $(n-1)$th Schr\"oder number~\cite{SS}, and that a permutation is separable if and only if it avoids the patterns $2413$ and $3142$~\cite{bo}. We shall generalize these results for $d$-permutations.

\subsection{Enumeration}
\label{sec:perm_enum}
In this section we generalize the result that the number of separable permutations of $[n]$ is the $(n-1)$th Schr\"oder number. Recall that the generating function of Schr\"oder numbers is $f = 1 + xf + xf^2$ \cite{oeis}.
\begin{sepenum}
\label{the:sepenum}
The generation function counting the number of separable $d$-permutations of $[n+1]$ satisfies
$f = 1 + xf + (2^{d-1}-1)xf^2$. The number $a_d(n)$ of $d$-permutations of $[n]$
satisfies the recursive formula:
$a_d(1)=1$, and for $n>1$
\begin{equation}\label{eq:enum_rec}
a_d(n) = 2^{d-1} \cdot \left(a_d(n-1) + \sum_{k=1}^{n-1} \frac{2^{d-1}-1}{2^{d-1}}a_d(k)a_d(n-k-1)\right),
\end{equation}
and for $n>1$ it is given by the formula
\begin{equation}\label{eq:enum_exp}
a_d(n) =
\frac{1}{n-1}
\sum_{k=0}^{n-2}
\binom{n-1}{k}\binom{n-1}{k+1}
(2^{d-1}-1)^k(2^{d-1})^{n-k-1}
.\end{equation}
\end{sepenum}

\noindent\textit{Proof.} Let $d$ be fixed.
For $n=1$ there is one $d$-permutation, which is separable.

Let $n>1$. Let $P$ be a separable $d$-permutation, and assume that its primary block structure is $(1 \ 2)$ in each row.
Consider its separation with the minimal possible $\ell$. If $\ell > 1$ then the left block is a separable $d$-permutation with primary block structure which has $(2 \ 1)$ at least in one row (because of the minimality of $\ell$); if $\ell=1$ then the primary block structure of the left block is not defined. The right block may be any separable $d$-permutation. Such a decomposition is unique, therefore, taking in account all $2^{d-1}$ possible primary block structures, we get
\[f=1+2^{d-1}xf\cdot\left(1+ \frac{2^{d-1}-1}{2^{d-1}}(f-1) \right)
= 1 + xf + (2^{d-1}-1)xf^2,\]
and the recursive formula (\ref{eq:enum_rec}) is clear from the same reasoning.

Ackerman \emph{et al.}~\cite{ack} obtained this recursive formula (with $d$ instead of $2^{d-1}$, and a(0)=1)
in their study of guillotine partitions. From this recursive formula they deduced an explicit formula which, after replacing $d$ by $2^{d-1}$ and $n$ by $n-1$, gives (\ref{eq:enum_exp}). We shall go into details on the connection between separable $d$-permutations and guillotine partitions in Section~\ref{sec:par_per}.
$\hfill\square$

\subsection{Characterization in terms of forbidden patterns}
\label{sec:perm_forbid}
Recall that a permutation is separable if and only if it avoids the patterns $2 4 1 3$, $3 1 4 2$. We generalize this result for separable $d$-permutations.

We shall use the following convention.
Let $P$ be a $d$-permutation represented by matrix. Take a restriction of $P$ to some $d'$ rows and $n'$ columns,
apply an order preserving relabeling on all the rows so that they will contain the numbers from $1$ to $n'$,
and exchange columns so that the first row will be $12\dots n'$. Denote the obtained matrix by $Q$. We say that $P$ contains $Q$ as a pattern.

In our discussion on separable permutations, we shall agree that the patterns are row-invariant, that is: A pattern $\pi$ is a $d'$-permutation, and any pattern obtained from $\pi$ by interchanging rows or columns is considered identical to $\pi$ (recall that we interchange columns in order to cause the first row be the natural order permutation).

\begin{forpat}
\label{the:forpat}
If a $d$-permutation $P$ contains any of the patterns
 \[
  \pi_1 = \left(
    \begin{array}{cccc}
      1 & 2 & 3 & 4\\ 
      2 & 4 & 1 & 3\\
    \end{array}
  \right), \ \
  \pi_2 = \left(
    \begin{array}{ccc}
      1 & 2 & 3 \\ 
      2 & 1 & 3 \\ 
      1 & 3 & 2 \\
    \end{array}
  \right), \ \
\pi_3 = \left(
    \begin{array}{ccc}
      1 & 2 & 3 \\ 
      2 & 3 & 1 \\ 
      3 & 1 & 2 \\
    \end{array}
  \right),\]
  then $P$ is non-separable.
\end{forpat}

\begin{figure}[ht]
$$\resizebox{45mm}{!}{\includegraphics{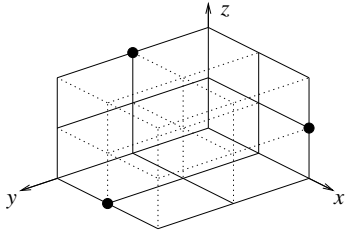}}$$
\caption{A non-separable $3$-permutation of $\{1, 2, 3\}$} \label{fig:nonsep3}
\end{figure}

This observation follows from the simple fact that any pattern in a separable $d$-permutation must be separable itself, and it is easy to check directly that $\pi_1$, $\pi_2$, $\pi_3$ are not separable. Observe that $\pi_1$ may be written in two forms:
\[\left(
     \begin{array}{cccc}
       1 & 2 & 3 & 4 \\
       2 & 4 & 1 & 3 \\
     \end{array}
   \right), \ \
  \left(
     \begin{array}{cccc}
       1 & 2 & 3 & 4 \\
       3 & 1 & 4 & 2 \\
     \end{array}
   \right);
   \]
$\pi_2$ may be written in six forms:
    \[
  \left(
    \begin{array}{ccc}
      1 & 2 & 3 \\
      2 & 1 & 3 \\
      1 & 3 & 2 \\
    \end{array}
  \right), \ \
  \left(
    \begin{array}{ccc}
      1 & 2 & 3 \\
      2 & 1 & 3 \\
      3 & 1 & 2 \\
    \end{array}
  \right), \ \
  \left(
    \begin{array}{ccc}
      1 & 2 & 3 \\
      2 & 3 & 1 \\
      1 & 3 & 2 \\
    \end{array}
  \right), \ \
  \left(
    \begin{array}{ccc}
      1 & 2 & 3 \\
      1 & 3 & 2 \\
      2 & 1 & 3 \\
    \end{array}
  \right), \ \
    \left(
    \begin{array}{ccc}
      1 & 2 & 3 \\
      1 & 3 & 2 \\
      2 & 3 & 1 \\
    \end{array}
  \right), \ \
  \left(
    \begin{array}{ccc}
      1 & 2 & 3 \\
      3 & 1 & 2 \\
      2 & 1 & 3 \\
    \end{array}
  \right);
  \]
$\pi_3$ may be written in two forms:
\[
  \left(
    \begin{array}{ccc}
      1 & 2 & 3 \\
      2 & 3 & 1 \\
      3 & 1 & 2 \\
    \end{array}
  \right), \ \
  \left(
    \begin{array}{ccc}
      1 & 2 & 3 \\
      3 & 1 & 2 \\
      2 & 3 & 1 \\
    \end{array}
  \right).
\]
The patterns $\pi_2$ and $\pi_3$, in all their forms, are all the patterns with $d=n=3$ which have a row with $2$ in the first position and a row with $2$ in the third position. Geometrically, these eight patterns (six representatives of $\pi_2$ and two representatives of $\pi_3$) are reflections and rotations of each other; Fig.~\ref{fig:nonsep3} presents one of them.

The next theorem is a characterization of separable $d$-permutations. It says that the patterns from Observation~\ref{the:forpat} may be taken as the only forbidden patterns.

\begin{sep3char}
A $d$-permutation $P$ of $[n]$ is separable if and only if it avoids the patterns
 \[
  \pi_1 = \left(
    \begin{array}{cccc}
      1 & 2 & 3 & 4\\ 
      2 & 4 & 1 & 3\\
    \end{array}
  \right), \ \
  \pi_2 = \left(
    \begin{array}{ccc}
      1 & 2 & 3 \\ 
      2 & 1 & 3 \\ 
      1 & 3 & 2 \\
    \end{array}
  \right), \ \
\pi_3 = \left(
    \begin{array}{ccc}
      1 & 2 & 3 \\ 
      2 & 3 & 1 \\ 
      3 & 1 & 2 \\
    \end{array}
  \right).\]
\end{sep3char}

\noindent\textit{Proof.}
The \emph{only if} direction is precisely Observation~\ref{the:forpat}, we shall prove the \emph{if} direction.

Let $d$ be fixed.

For $n \geq 3$ the proof is by induction on $n$. For $n = 3$, it is easy to see that if $2$ does not appear in the first column or does not appear in the last column (and this precisely means that $P$ avoids $\pi_2$ and $\pi_3$) then $P$ is separable with, respectively, $\ell=1$ or $\ell=2$.
Note that according to the formula from Observation~\ref{the:sepenum},
the number of separable $d$-permutations of
$\{1, 2, 3\}$ is $2^{d-1}(2^d-1)$. It is easy to see that this is precisely the number of ways
to write a $d$-permutation of $\{1, 2, 3\}$ with $2$s only in two first columns or only in two last columns.


Let $P$ is a $d$-permutation of $[n+1]$. We shall show that either $P$ is separable, or it contains one of the patterns $\pi_1$, $\pi_2$, $\pi_3$. This suffices in view of Observation~\ref{the:forpat}.

Apply on all the rows of $P$ the order preserving relabeling so that first $n$ columns contain members of $[n]$, and the numbers in $P^{(n+1)}$ (the $(n+1)$th column of $P$) belong to $\{\frac{1}{2}, 1\frac{1}{2}, \dots, n+\frac{1}{2}\}$. \emph{Example:}
\[
\left(
  \begin{array}{ccccc}
    1 & 2 & 3 & 4 & 5 \\
    4 & 5 & 1 & 2 & 3 \\
    5 & 2 & 4 & 3 & 1 \\
  \end{array}
\right)
\rightarrow
\left(
  \begin{array}{ccccc}
\vspace{1mm}    1 & 2 & 3 & 4 & 4\frac{1}{2} \\
\vspace{1mm}    3 & 4 & 1 & 2 & 2\frac{1}{2} \\
    4 & 1 & 3 & 2 & \frac{1}{2} \\
  \end{array}
\right).
\]

We shall keep calling this object $P$. The matrix formed by its first $n$ columns is a $d$-permutation of $[n]$ which we denote by $P'$. If $P'$ is not separable, then it contains one of the forbidden patterns by induction hypothesis, therefore $P$ also contains it.

Thus we suppose from now on that $P'$ is separable. Let $\ell$ be the \emph{minimal} number, $1 \leq \ell < n$, so that $P'$ is separated between $\ell$ and $\ell+1$. We assume that the primary block structure of $P'$ is $(1 \ 2)$ in each row: there is no loss of generality because otherwise we can relabel the members of any row according to $j \leftrightarrow (n-j+1)$; it is easy to see that separability of a $d$-permutation, and avoiding the patterns $\pi_1$, $\pi_2$, $\pi_3$ are invariant under this transformation (which is, geometrically, reflection with respect to the direction of an axis).

Consider $P^{(n+1)}$. If all the numbers in $P^{(n+1)}$ belong to $\{\ell+\frac{1}{2}, \dots, n+\frac{1}{2}\}$, then $P$ may be separated between $\ell$ and $\ell+1$. If all the numbers in $P^{(n+1)}$ belong to $\{\frac{1}{2}, n+\frac{1}{2}\}$, then $P$ may be separated between $n$ and $n+1$. Thus we assume from now on that there is a member of $P^{(n+1)}$ that belongs to $\{\frac{1}{2}, \dots, \ell-\frac{1}{2}\}$, and not all of them are $\frac{1}{2}$ or $n+\frac{1}{2}$.

Suppose that one of the members in $P^{(n+1)}$ is $\frac{1}{2}$ -- assume that this happens in the row $p_2$.
According to our assumption, $P^{(n+1)}$ contains a number that belongs to $\{1+\frac{1}{2}, \dots, n-\frac{1}{2}\}$ -- suppose that this happens in the row $p_3$. Restrict $P$ to the rows $p_1, p_2, p_3$ and the following three columns: (1) the column which contains $1$ in $p_3$, (2) the column which contains $n$ in $p_3$, (3) $P^{(n+1)}$. This restriction has the following form -- we write either an exact number, or the interval to which it belongs:

\[ \left(
     \begin{array}{ccc}
       [1\dots \ell] & [\ell+1\dots n] & n+\frac{1}{2} \\

       [1 \dots \ell] & [\ell+1 \dots n] & \frac{1}{2} \\

       1 & n & [1+\frac{1}{2} \dots n-\frac{1}{2}]
     \end{array}
   \right),
\]
or, after relabeling,
\[
\left(
  \begin{array}{ccc}
    1 & 2 & 3 \\
    2 & 3 & 1 \\
    1 & 3 & 2 \\
  \end{array}
\right),
\]
which is (a form of) $\pi_2$.

In particular, now the case $\ell=1$ is settled, since assumed that $P^{(n+1)}$ contains a number from the interval $\{\frac{1}{2}, \dots, \ell-\frac{1}{2}\}$. And we assume from now on that $P^{(n+1)}$ does not contain $\frac{1}{2}$ and that $\ell>1$.

Recall that not all the numbers in $P^{(n+1)}$ are in $\{\ell+\frac{1}{2}, \dots, n+\frac{1}{2}\}$. Therefore it contains a number from $\{1+\frac{1}{2}, \dots, \ell-\frac{1}{2}\}$ -- suppose that this happens in the row $p_2$.

Suppose $P'(L)$ be separated between $k$ and $k+1$,
where $1\leq k < \ell$ (recall that $\ell>1$). Consider the primary block structure of $P'(L)$. Suppose first that it is $(2 \ 1)$ in the row $p_2$. Take the restriction of $P$ to the rows $p_1$ and $p_2$ and the following four columns: (1) the column which has $1$ in $p_2$, (2) the column which has $\ell$ in $p_2$, (3) the column which has $n$ in $p_2$, (4) $P^{(n+1)}$. This restriction is
\[\left(
    \begin{array}{cccc}
      [1 \dots k] & [k+1 \dots \ell] & [\ell+1 \dots n] & n+\frac{1}{2} \\
      \ell & 1 & n & [1+\frac{1}{2} \dots \ell-\frac{1}{2}] \\
    \end{array}
  \right),
\]
or, after relabeling,
\[\left(
  \begin{array}{cccc}
    1 & 2 & 3 & 4 \\
    3 & 1 & 4 & 2 \\
  \end{array}
\right),\]
which is (a form of) $\pi_1$.

We assume from now on that the primary block structure of $P'(L)$ in the row $p_2$ is $(1 \ 2)$. However it cannot be $(1 \ 2)$ in all the rows, because the minimality in the choice of $\ell$. Therefore there is a row (say $p_3$), such that the primary block structure of $P'(L)$ in $p_3$ is $(2 \ 1)$. Consider $P_{3, n+1}$. If it belongs to $\{1+\frac{1}{2}, \dots, \ell-\frac{1}{2}\}$, we obtain a pattern $\pi_{1}$ as just discussed. Therefore we assume that it belongs to $\{\ell+\frac{1}{2}, \dots, n+\frac{1}{2}\}$.

We have here two cases. If $P_{3, n+1} > \ell+\frac{1}{2}$: take $P$ restricted to the rows $p_2$ and $p_3$ and to the columns: (1) the column which has $1$ in $p_2$, (2) the column which has $\ell$ in $p_2$, (3) the column which has $\ell+1$ in $p_3$, (4) $P^{(n+1)}$. This restriction is
\[\left(
    \begin{array}{cccc}
      1 & \ell & [\ell+1 \dots n] & [1+\frac{1}{2} \dots \ell-\frac{1}{2}] \\

      [\ell-k+1 \dots \ell] & [1 \dots \ell-k] & \ell+1 & [\ell+1+\frac{1}{2} \dots n+\frac{1}{2}] \\
    \end{array}
  \right),
\]
or, after relabeling,
\[\left(
  \begin{array}{cccc}
    1 & 3 & 4 & 2 \\
    2 & 1 & 3 & 4 \\
  \end{array}
\right),\]
which is (a form of) $\pi_1$.

In the second case, $P_{3, n+1} = \ell+\frac{1}{2}$: take $P$ restricted to the rows $p_1$, $p_2$ and $p_3$ and to the columns: (1) the column which has $1$ in $p_2$, (2) the column which has $\ell$ in $p_2$, (3) $P^{(n+1)}$. This restriction is
\[\left(
    \begin{array}{ccc}
     [1 \dots \ell-1] & [\ell \dots n] & n+\frac{1}{2} \\

      1 & \ell & [1+\frac{1}{2} \dots \ell-\frac{1}{2}] \\

      [\ell-k+1 \dots \ell] & [1 \dots \ell-k] & [\ell+1+\frac{1}{2} \dots n+\frac{1}{2}] \\
    \end{array}
  \right),
\]
or, after relabeling,
\[\left(
  \begin{array}{ccc}
    1 & 2 & 3  \\
    1 & 3 & 2  \\
    2 & 1 & 3  \\
  \end{array}
\right),\]
which is (a form of) $\pi_2$. $\hfill\square$

\section{Guillotine partitions}

\subsection{Introduction}
\label{sec:guil_intro}
For the sake of completeness, we remind the definition of guillotine partition. The next paragraph, containing this definition, is taken from~\cite{ack} almost verbatim, with only a slight change in notation:
\begin{quote}
Let $B$ be an axis-parallel $d$-dimensional box in $\mathbb{R}^d$. A
\emph{partition} of $B$ is a set $S$ of $k >0$ interior-disjoint
axis-parallel boxes $b_1, b_2, \dots, b_k$ whose union equals $B$. A
partition $S = \{b_1, b_2, \dots, b_k\}$ of $B$ is a \emph{guillotine partition}
if $k = 1$ or there are a hyperplane $h$ and two disjoint non-empty
subsets $S^-, S^+ \subset S$ such that:
\begin{enumerate}
  \item $h$ splits $B$ into two interior-disjoint boxes $B^-$ and $B^+$;
  \item $S^-$ is a guillotine partition of $B^-$;
  \item $S^+$ is a guillotine partition of $B^+$.
\end{enumerate}
\end{quote}

In this definition, the hyperplane $h$ is orthogonal to some
axis $x_i$. It is assumed that the interior of $B^-$ is below $h$
and the interior of $B^+$ is above $h$, with respect to $x_i$.

\begin{figure}[ht]
$$\resizebox{150mm}{!}{\includegraphics{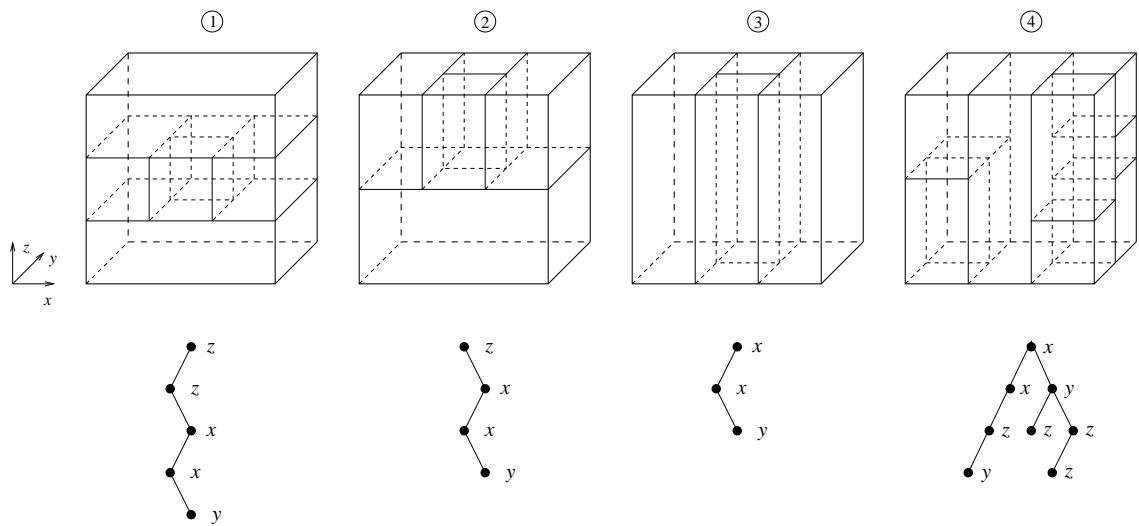}}$$
\caption{Guillotine partitions of a $3$-box and the corresponding trees}
\label{fig:bad_graphs}
\end{figure}

Ackerman \emph{et al.}~\cite{ack}
enumerated structurally different guillotine
partitions of $B$ by $n$ hyperplanes and established a
bijection between the set of such partitions and the set of binary trees with $n$ vertices, each vertex colored by a color belonging to the set $\{1, 2, \dots, d\}$, with the restriction: if a vertex
$v$ is a right child of $u$, then these vertices have different
colors.
The binary tree corresponding to a guillotine partition $S$ is constructed recursively
as follows. An empty tree corresponds to the trivial partition $S=\{B\}$.
Otherwise, consider a hyperplane $h$ that splits $B$ into two subboxes as in the definition.
If there are several such hyperplanes, they must be
orthogonal to the same axis $x_i$; in this case choose $h$ to be the
highest among them (with respect to $x_i$). The root of the
corresponding tree is then colored by $i$; the left branch of the
root is the tree that corresponds to the partition of $B^-$ and the
right branch of the root is the tree that corresponds to the
partition of $B^+$. The choice of $h$ implies that no vertex has a
right child with the same color, and it is easily proved
recursively that this correspondence is indeed a bijection. See Fig.~\ref{fig:bad_graphs} for some examples.

The main enumeration results in \cite{ack} are an exact formula for a number
of structurally different guillotine partitions of $B$ (see
\cite[A103209]{oeis}), and its asymptotic behavior.

We construct a bijection between
between the set guillotine partitions of a $2^{d-1}$-dimensional box by $n$ cuts and separable $d$-permutations of $[n+1]$.
In addition, we consider structurally different guillotine partitions under several natural restrictions, enumerate them by means of generating function or explicit formula,
find some connections with other combinatorial structures,
and use it for enumeration of permutations under certain restrictions.

\vspace{2mm}

\noindent\textit{Notation and convention.} ``The number of guillotine partitions'' stands for ``the number of structurally different guillotine partitions''.

The $d$-dimensional box (or \emph{$d$-box}, for short) being partitioned will be denoted by $B$. A \emph{subbox} of $B$ is a subset of $B$ obtained at some recursive stage of constructing a guillotine partition.

Let $h$ be a hyperplane which splits a subbox of $B$ into two smaller subboxes, as in the definition of guillotine partition. A \emph{cut} is the intersection of such a hyperplane with appropriate $(d-1)$-dimensional faces of these subboxes.

In a binary tree, a \emph{left descendant of a vertex $x$} is either the left
child of $x$ or a descendant of the left child of $x$.
A \emph{right descendant} is defined similarly.

Partitions are identified with trees that correspond to them under the bijection described above. Therefore consider cuts as vertices of the tree and we shall occasionally use expressions like ``right (left) child (descendant)'' for cuts in $B$.

If $h$ is a cut and we say \emph{the higher} (or \emph{lower}) \emph{half-space bounded by $h$}, this means higher (or lower) with respect to the axis to which $h$ is orthogonal.

The \emph{principal cut} of $B$ is a cut which splits it into two parts. It was noted above that all principal cuts of $B$ are parallel. Note that the \emph{highest principal cut} corresponds to the root of the tree.

On figures, the letter right to a vertex denotes its color, the letter left to it denotes its label.

In each section, $f$ denotes the generating function in the case discussed in this section. For fixed $d$, it is also denoted by $f_d$. The coefficient of $x^n$ in $f_d$ -- that is, the number of guillotine partitions of a $d$-box by $n$ cuts, with the relevant restriction -- will be denoted by $a_d(n)$.

\subsection{Schr\"{o}der paths}
For $d=2$ the counting sequence of guillotine partitions is the sequence of
Schr\"{o}der numbers
(this was also found by Yao \emph{et al.} \cite{yao}). We remind how this follows from considering the generating function, and construct a bijection between guillotine partitions and an appropriate generalization of Schr\"{o}der paths -- for general $d$.

Let $f$ be the generating function for the number of guillotine
partitions of a $d$-box.
Since all possible guillotine partitions may take place in $B^-$, and only those with principal cut in a different direction -- in $B^+$,
we have $f = 1+dxf^-f^+$, where $f^-=f$
and $f^+ = 1+\frac{d-1}{d}(f-1)$. It follows $f=1+xf+(d-1)xf^2$. For
$d=2$ this is the generating function of Schr\"{o}der
numbers. The general case may be interpreted as
follows:

\begin{schr}
\label{the:schr} There is a bijection between the set $S_{d, n}$ of
Schr\"{o}der paths of length $2n$ with up-steps colored by $\{1,
\dots, d-1\}$ and the set $T_{d, n}$ of binary trees with $n$
vertices colored by $\{0, 1, \dots, d-1\}$ avoiding a vertex and its
right child of the same color.
\end{schr}

In a Schr\"oder path, letters $\y U$, $\y D$, and $\y L$ denote up-steps, down-steps and level-steps, respectively.

\vspace{2mm}

\noindent \textit{Proof of Observation~\ref{the:schr}.} We construct a
bijection $\varphi$ from $\cup_{d\geq 1, n\geq 0 } S_{d,n}$ to
$\cup_{d\geq 1, n\geq 0 } T_{d,n}$ as follows.

For $n=0$: the empty tree corresponds to the empty Schr\"{o}der path
($\varphi(\emptyset)=\emptyset$).

For $n \geq 1$: Each Schr\"{o}der path $P$ may be decomposed in
precisely one of the three following ways: $P=\y{L}Q$,
$P=\y{U}Q\y{D}$ (where $Q$ is a Schr\"{o}der path of length
$2(n-1)$), or $P=\y{U}Q\y{D}R$ (where $Q$ and $R$ are Schr\"{o}der
paths of total length $2(n-1)$ and $R$ is non-empty). Define $\varphi(P)$ as follows:
\begin{itemize}
  \item If $P=\y{L}Q$: the left branch of $\varphi(P)$ is $\varphi(Q)$, and the root of $\varphi(P)$ is colored by $0$.
  \item If $P=\y{U}Q\y{D}$ and $\y{U}$ is colored by $a$: the left branch of $\varphi(P)$ is $\varphi(Q)$, and the root of $\varphi(P)$ is colored by $a$.
  \item If $P=\y{U}Q\y{D}R$ and $\y{U}$ is colored by $a$: the left branch of $\varphi(P)$ is $\varphi{Q}$, the left branch of $\varphi(P)$ is $\varphi(R)$ and the root of $\varphi(P)$ is colored by $(a+b) (\y{mod\ }d)$ where $b$ is the color of the root of $\varphi(Q)$.
\end{itemize}

No vertex and its right child
colored by the same color in the tree $\varphi(P)$, since in the third case $a \not = 0$.

\begin{figure}[ht]
$$\resizebox{120mm}{!}{\includegraphics{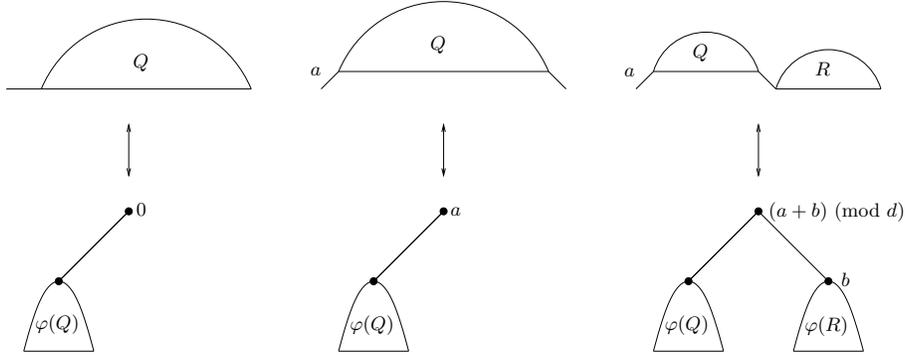}}$$
\caption{Proof of Observation~\ref{the:schr}}
\label{fig:schr}
\end{figure}

It is easy to see that $\varphi$ is
bijective, and that the image of its restriction to $S_{d, n}$ is
$T_{d, n}$. See Fig.~\ref{fig:schr} for illustration.

\subsection{Guillotine partitions and separable $d$-permutations}
\label{sec:par_per}
As we already mentioned in the proof of Observation~\ref{the:sepenum}, Ackerman \emph{et al.} \cite{ack} found a recursive formula for the number of guillotine partitions of a $d$-box by $n$ cuts, which, after replacing $d$ by $2^{d-1}$ and $n$ by $n-1$ gives the formula for the number of $d$-permutations of $[n+1]$. In the next theorem we construct a bijection between these structures.

\begin{permd}
\label{the:permd}
There is a bijection between the set of guillotine partitions of $2^{d-1}$-dimensional box by $n$ cuts and the set of separable $d$-permutations of $[n+1]$.
\end{permd}

\noindent\textit{Proof.} Consider any bijective correspondence between the axes of $\mathbb{R}^{2^{d-1}}$ and
all possible primary block structures of $d$-permutations
(geometrically -- the pairs of orthants in $\mathbb{R}^d$) (a canonical correspondence will use the binary representation; however, we won't use the particular form of the correspondence).

For $n=0$, the trivial partition corresponds to the only $d$-permutation of $\{1\}$.

Let $P$ be a partition of $B$ with $n>0$ cuts. Let $h$ be its highest principal cut, and assume that it is orthogonal to $x_i$-axis. Construct the (unique) $d$-permutation which has the $d$-permutation corresponding to $B^-$ as the left block, the $d$-permutation corresponding to $B^+$ as the right block, and has the primary block structure corresponding to $x_i$-axis. (Geometrically, we put the graphs of the $d$-permutations corresponding to $B^-$ and $B^+$ into an appropriate orthants.) This $d$-permutation corresponds to $P$. It is easy to see that this correspondence is a bijection. $\hfill\square$

\vspace{2mm}

\begin{figure}[ht]
$$\resizebox{110mm}{!}{\includegraphics{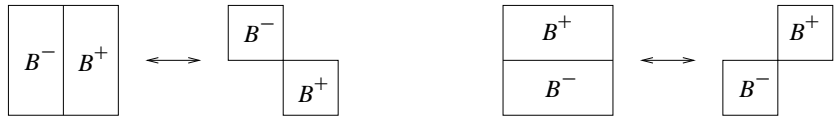}}$$
\caption{Constructing a permutation from a planar guillotine sequence, a recursive step}
\label{fig:perm1}
\end{figure}

\begin{figure}[ht]
$$\resizebox{150mm}{!}{\includegraphics{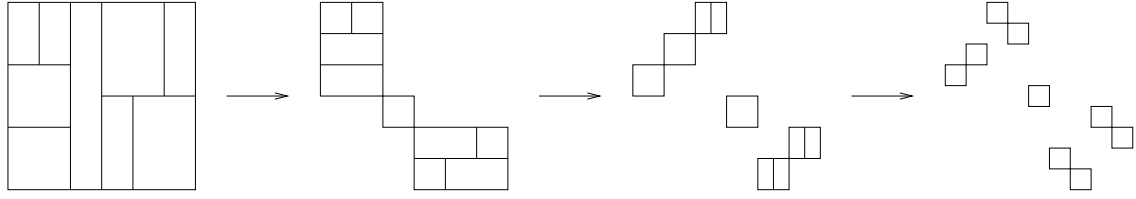}}$$
\caption{Constructing a permutation from a planar guillotine sequence, an example}
\label{fig:perm3}
\end{figure}

The special case $d=2$ has an especially clear visualization.
Let horizontal principal cuts correspond to the primary block structure
$\left(
\begin{array}{cc}
1 & 2 \\
1 & 2 \\
\end{array}
\right)$, and vertical cuts to
$\left(
\begin{array}{cc}
1 & 2 \\
2 & 1 \\
\end{array}
\right)$.
This means: Given a partition, if a highest principal cut is vertical we slide the parts along the cut until the left part is above the right part, and if a highest principal cut is horizontal we slide the parts along the cut until the upper part is right to the lower part, see Fig.~\ref{fig:perm1}.
Continuing this recursively, we obtain at some stage boxes with trivial partitions. Replacing them with points, we obtain the graph of the corresponding permutation. See Fig.~\ref{fig:perm3} for an example.

In fact, our bijection in this case ($d=2$) is a version of a special case of a bijection found by
Ackerman \emph{et al.} \cite{ack3, ack2}.
They found a bijection between all planar rectangular partitions of a square by $n$ cuts and Baxter permutations of $[n+1]$. A restriction of this bijection is a correspondence between guillotine partitions and separable permutations, which is essentially equal to our correspondence.
More precisely, the permutation corresponding to a guillotine partition under our bijection and the permutation obtained from the same partition by applying FP2BR algorithm from \cite{ack3}, are obtained from each other by relabeling $j \leftrightarrow (n-j+1)$.

Recall that Baxter permutations are those avoiding $25\bar314$ and $41\bar352$ and thus separable permutations form their subfamily.
Thus the $2$-dimensional case of Theorem~\ref{the:permd} may be interpreted as follows: While all partitions correspond to Baxter permutations which allow patterns $2413$ and $3142$ only on some condition, guillotine partitions correspond precisely to those avoiding these patterns. Thus we have here the following problem for a future research: \emph{to describe a class of $d$-permutations which would generalize Baxter permutations in the sense of being in a natural bijection with $2^{d-1}$-dimensional guillotine partitions}.

Let us return to Theorem~\ref{the:permd}.
It describes a correspondence between partitions of a $q$-dimensional box and multidimensional permutations when $q$ is a power of $2$. Can we find a similar correspondence when $q$ is not a power of $2$? Yes: let us take $d$ so that $2^{d-1} \geq q$ and to use only $q$, among $2^{d-1}$, primary block structures in forming separated $d$-permutations. The same reasoning as in the proof of Theorem~\ref{the:permd}
gives then the following result:

\begin{permq}
\label{the:permq}
For $q \leq 2^{d-1}$,
there is a bijection between
the set of guillotine partitions of $q$-dimensional box by $n$ cuts and
the set of separable $d$-permutations of $[n+1]$
avoiding any $2^{d-1}-q$ (among $2^{d-1}$) $d$-permutations of $\{1, 2\}$.
\end{permq}

\section{Restricted guillotine partitions}
\label{sec:restr}
In this section we consider guillotine partitions with certain natural restrictions.
We find their generating functions, in some cases -- explicit formulae, and show connections with other combinatorial structures. In planar case, this will also help us to obtain several results on patterns in permutations.

\subsection{Boundary guillotine partitions}
\label{sec:bound}

Consider the following subfamily of guillotine partitions.
Let $B$ be a $d$-box. Each cut in a guillotine partition
of $B$ is a $(d-1)$-dimensional box which has $d-1$ pairs of
opposite $(d-2)$-dimensional faces. We require that for all cuts, in each such
pair at least one face belongs to the boundary of $B$. A guillotine
partition that satisfies this condition will be called a \emph{boundary
guillotine partition.}
\emph{Example:} Among four guillotine partitions on Fig.~\ref{fig:bad_graphs},
only (4) is a boundary guillotine partition.

\begin{guil2}

\label{the:guil2}
For fixed $d$, the generating function counting the number of boundary guillotine partitions of a $d$-dimensional box obtained by $n$ cuts is
\[f = 1 + \frac{dx(1-x)}{(1-x)^2}\left( 1 + \frac{(d-1)x(1-x)}{(1-2x)^2} \left(\dots \left( 1 + \frac{2x(1-x)}{(1-(d-1)x)^2}\left( 1 + \frac{x(1-x)}{(1-dx)^2} \right)^2\right)^2 \dots \right)^2 \right)^2.\]
\end{guil2}

\noindent \textit{Proof.}
Let $B$ be a $d$-dimensional box with a boundary guillotine partition. For each $i = 1, \dots, d-1$
we define inductively subboxes $B(\underbrace{--\dots-}_{i})$,
and $B(\underbrace{--\dots-}_{i-1}+)$ $B^*(\underbrace{--\dots-}_{i})$ as follows.
Let $h_1$ be the highest principal cut of $B$. Assume without loss of generality that $h_1$ is orthogonal to $x_1$-axis. Denote by $B(-)$ the part
of $B$ below $h_1$, and by $B(+)$ the part of $B$ above $h_1$ (thus they are just $B^-$ and $B^+$ from the definition of guillotine partition). $B(-)$ may
have several principal
cuts parallel to $h_1$. Denote by $B^*(-)$ the part of $B$ below the lowest of them.

Let $1 < i \leq d$.
Consider $B^*(\underbrace{--\dots-}_{i-1})$.
By induction, it is bounded form above by cuts orthogonal to $x_1$, $x_2$, \dots, $x_{i-1}$, and its highest principal cut $h_{i}$ is orthogonal to a ``new'' axis -- $x_{i}$ without loss of generality.
The cut $h_{i}$ splits $B^*(\underbrace{--\dots-}_{i-1})$ into two subboxes:
$B(\underbrace{--\dots-}_{i})$ below $h_{i}$ , and $B(\underbrace{--\dots-}_{i-1}+)$ above $h_{i}$.

\begin{figure}[ht]
$$\resizebox{160mm}{!}{\includegraphics{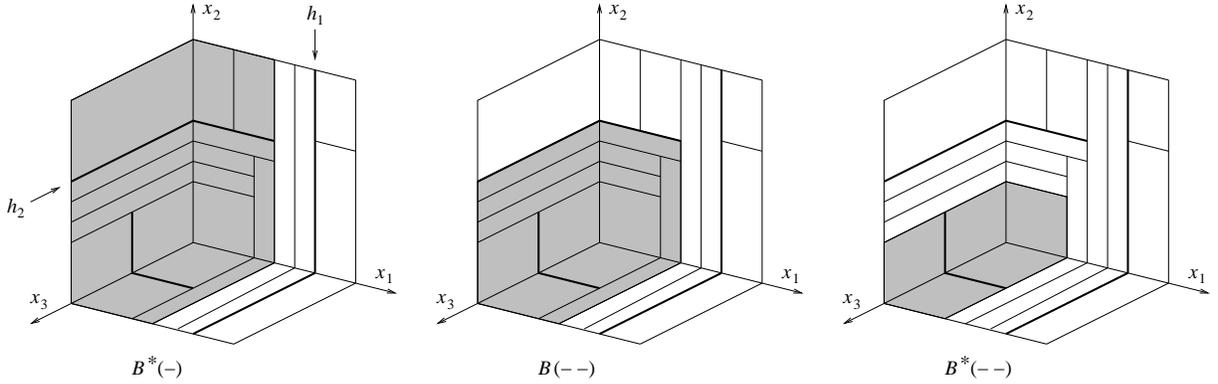}}$$
\caption{Distinguished subboxes of $B$ in the proof of Theorem~\ref{the:guil2}} \label{fig:b}
\end{figure}

If $i<d$: Consider $B(\underbrace{--\dots-}_{i})$. Starting with $h_{i}$, pass to its left child (if it exists), then to its left child, and so on until the first occurrence of a cut (denote it by $h_{i+1}$) which is not orthogonal to $x_1$, $x_2$, \dots, $x_{i}$. Denote by $B^*(\underbrace{--\dots-}_{i})$ the subbox whose highest principal cut is $h_{i+1}$.
See Fig.~\ref{fig:b}.

Denote by $f(\underbrace{--\dots-}_i)$, $f(\underbrace{--\dots-}_{i-1}+)$
the generating functions for number of boundary guillotine partitions of $B(\underbrace{--\dots-}_i)$, $B(\underbrace{--\dots-}_{i-1}+)$,
respectively.
We have $f = 1 + dxf(-)f(+)$ and $f(-)=\frac{1}{1-x}f(+)$, 
and therefore $f = 1 + dx(1-x)\cdot(f(-))^2$.

For $1\leq i < d$ we have
\[f(\underbrace{--\dots -}_{i}) = \frac{1}{1-ix} \left( 1+(d-i)xf(\underbrace{--\dots -}_{i+1})f(\underbrace{--\dots -}_{i}+) \right).\]
The factor $\frac{1}{1-ix}$ is due to the fact that before we reach $B^*(\underbrace{--\dots-}_{i})$ we (possibly) have cuts orthogonal to $x_j$, $j \leq i$, that may appear in any possible order which determines their structure completely. The factor $(d-i)$ is due to the fact that the cut in ``new'' direction may be orthogonal to any $x_j$, $i < j \leq d$.

In addition, $f(\underbrace{--\dots -}_{i+1}) = \frac{1}{1-x}f(\underbrace{--\dots -}_{i}+)$ because of possible several principal cuts, and thus

\[f(\underbrace{--\dots -}_{i}) = \frac{1}{1-ix} \left( 1+(d-i)x(1-x)\cdot\left(f(\underbrace{--\dots -}_{i+1})\right)^2 \right).\]

Finally, $f(\underbrace{--\dots -}_{d}) = \frac{1}{1-dx}$ because there is no more unused direction, and this implies $f$ as stated in the theorem. $\hfill\square$

\vspace{2mm}

Table~\ref{tab:boundary} presents the number of boundary guillotine partitions for $2 \leq d \leq 5$, $0 \leq n \leq 15$.

For $d=2$, we have
\[f = 1+\frac{2x(1-3x+3x^2)^2}{(1-x)(1-2x)^4}.\]
Finding the coefficients of this function, we get the explicit formula: the number of boundary partitions of a square obtained by $n$ cuts is
\[a_2(n) = 2+\frac{(n-1)(n^2+n+42)}{3}2^{n-4}\]
when $n\geq 1$, and $a_2(0)=1$.

\begin{table}[ht]
\begin{center}
\begin{tabular}{|c|c|c|c|c|}
 \hline
 $n$ & $a_2(n)$ & $a_3(n)$ & $a_4(n)$ & $a_5(n)$\\ \hline
 0 & 1     &1          & 1         &1\\
 1 & 2     &3          & 4         &5\\
 2 & 6     &15         &28         &45 \\
 3 & 20    &87         &232        &485  \\
 4 & 64    &507        &1984       &5485  \\
 5 & 194   &2859       &16804      &62405  \\
 6 & 562   &15495      &139012     &702445  \\
 7 & 1570  &80943      &1119172    &7770085  \\
 8 & 4258  &409539     &8776036    &84292525  \\
 9 & 11266 &2015907    &67190308   & 897101125  \\
 10& 29186 &9687855    &503591332  &  9379187885\\
 11& 74242 &45574791   &3703703716 &  96487985125\\
 12&185858 &210305739  &26779859332&  978249364205\\
 13&458754 &953479899  &190652265220& 9787794765765\\
 14&1118210&4252898199 &1337960522980&96752629782125\\
 15&2695170&18683733663&9264356481220&945738292868325\\
\hline
\end{tabular}
\vspace{2mm}
\caption{Values of $a_d(n)$ for $2 \leq d \leq 5$, $0 \leq n \leq 15$}
\label{tab:boundary}
\end{center}
\end{table}

We can also estimate the asymptotic behavior of $f_d$. Recall that sequences $a_n$ and $b_n$ are \emph{asymptotically equivalent}
as $n\rightarrow\infty$ if
$\lim_{n\rightarrow\infty}\frac{a_n}{b_n}=1$; this is denoted by $a_n\sim
b_n$. Theorem~\ref{the:guil2} implies the following result:

\begin{guil3}
For all $d\geq2$,
\[[x^n]f_d(x)\sim\left(\frac{d-1}{d}\right)^{2^d-1}\frac{n^{2^d-1}}{(2^d-1)!}d^n.\]
\end{guil3}

\noindent\textit{Proof.} Define
$$g_d(x)=
\left(\frac{dx(1-x)}{(1-x)^2}\right)^{2^0}
\left(\frac{(d-1)x(1-x)}{(1-2x)^2}\right)^{2^1}
\cdots
\left(\frac{2x(1-x)}{(1-(d-1)x)^2}\right)^{2^{d-2}}
\left(\frac{x(1-x)}{(1-dx)^2}\right)^{2^{d-1}},$$ which is
equivalent to
$$g_d(x)=x^{2^d-1}(1-x)^{2^d-1}
\frac{\prod_{j=0}^{d-1}(d-j)^{2^j}}{\prod_{j=1}^d(1-jx)^{2^j}}.$$
It
is not hard to see that $[x^n]f_d(x)\sim [x^n]g_d(x)$, and the
smallest positive pole of the function $g_d(x)$ is $x^*=\frac{1}{d}$
of order $2^{d}$. Hence,
$$[x^n]f_d(x)\sim c_d\frac{n^{2^d-1}}{(2^d-1)!}d^n,$$
where

\[c_d=\lim_{x\rightarrow x^*}(x-x^*)^{2^d}g_d(x)
=(x^*)^{2^d-1}(1-x^*)^{2^d-1}\frac{\prod_{j=0}^{d-1}(d-j)^{2^j}}{\prod_{j=1}^{d-1}(1-jx^*)^{2^j}}
=\frac{(d-1)^{2^d-1}}{d^{2^{d+1}-2}}\frac{\prod_{j=0}^{d-1}(d-j)^{2^j}}{\prod_{j=1}^{d-1}\frac{(d-j)^{2^j}}{d^{2^j}}}=\]
\[=\frac{(d-1)^{2^d-1}}{d^{2^{d+1}-1-2^0-2^1-2^2\cdots-2^{d-1}}}\frac{\prod_{j=0}^{d-1}(d-j)^{2^j}}{\prod_{j=1}^{d-1}(d-j)^{2^j}}
=\frac{d(d-1)^{2^d-1}}{d^{2^d}}
=\frac{(d-1)^{2^d-1}}{d^{2^d-1}},\]
as claimed.$\hfill\square$

\vspace{2mm}
Recall the bijection between guillotine partitions and binary trees (see Section~\ref{sec:guil_intro}). We characterize the trees corresponding to boundary partitions under this bijection.

\begin{guil}
\label{the:forb} There is a bijection between boundary guillotine
partitions of a $d$-box by $n$ hyperplane cuts and binary trees with
$n$ vertices colored by $\{1, 2, \dots, d\}$ that satisfy the
following:
\begin{itemize}
  \item If $v$ is a right child of $u$, then they have different colors.
  \item If $v$ is a right descendant of $u$, $w$ is the left child of $v$,
$v$ and $u$ having the same color, then the color of $w$ is different from the color of $u$ and $v$.
  \item If $v$ is a left descendant of $u$, $w$ is the right child of $v$,
$v$ and $u$ having the same color, then the color of $w$ is different from the color of $u$ and $v$.

\end{itemize}
\end{guil}

\begin{figure}[ht]
$$\resizebox{60mm}{!}{\includegraphics{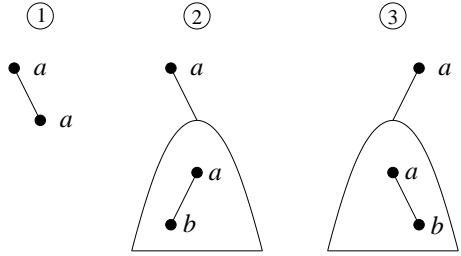}}$$
\caption{Forbidden types of trees for boundary partitions} \label{fig:forb}
\end{figure}

Figure~\ref{fig:forb} presents the three forbidden types of
subtrees.

\noindent \textit{Proof of Observation~\ref{the:forb}.} The subtree of
type (1) is forbidden in binary trees corresponding to all
guillotine partition, and this is the only forbidden subtree in the
general case~\cite[Observation 2]{ack}.

Suppose that the tree corresponding to a partition has a subtree of
type (2). The cut $u$ partitions a subbox of $B$ into two
subboxes.
The cut $v$ lies
above $u$. Therefore the cut $w$
is bounded by the cut
$v$ from above and by the cut $u$ from below (with respect to the $x_a$-axis). Thus both
$(d-2)$-dimensional faces of $w$ orthogonal to $x_a$ do not belong
to the boundary of $B$: they belong to $u$ and $v$.
A similar reasoning proves that the subtrees of type (3) are
forbidden.

Suppose now that a guillotine partition $P$ is not a boundary partition.
Assume without loss of generality that $P$ is a minimal not boundary
guillotine partition. That is: $P$ is not a boundary partition, but it is obtained
by joining two boundary partitions -- say, those of boxes $B^-$ and $B^+$, along a cut
$u$ orthogonal to the $x_a$-axis. There is a cut $w$ in
one of the boxes -- assume that $w$
is in $B^+$ and it is orthogonal to $x_b$-axis -- which meets
$u$ (from above) and
another cut $v$ orthogonal to $h_a$ (from below).
Then $w$ is the left
child of $v$, and $v$ is a right descendant of $u$. Besides, $u$ and $v$
have the color $a$, and $w$ has the color $b$, thus we obtain a
subtree of type (2). Similarly, assuming $w$ is in $B^-$ we obtain
a subtree of type (3). $\hfill\square$

\begin{figure}[ht]
$$\resizebox{80mm}{!}{\includegraphics{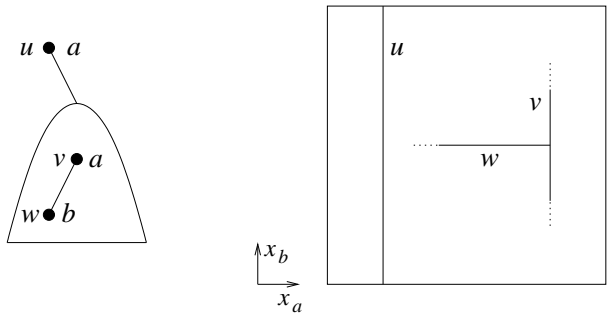}}$$
\caption{Illustration to the proof of Observation~\ref{the:forb}} \label{fig:uvw1}
\end{figure}

\vspace{2mm} Figure~\ref{fig:uvw1} illustrates the proof of
Observation~\ref{the:forb}.

\subsection{Alternating guillotine partitions}
\label{sec:alter}

In this section we consider another restriction on guillotine
partitions. Fix $m\geq 2$. A guillotine partition is \emph{$m$-alternating}
if each subbox has at most $(m-1)$ principal cuts.

\begin{alt0}
\label{the:alt0}
For fixed $d$, the generating function $f$ counting the number of $m$-alternating guillotine partitions of a $d$-dimensional box obtained by $n$ cuts
satisfies
\[f = \left(1 + \frac{d-1}{d}(f-1)\right)
\left(
\frac
{1 - x^m \left( 1+\frac{d-1}{d}(f-1) \right)^m}
{1 - x \left( 1+\frac{d-1}{d}(f-1) \right)}
\right) .\]
\end{alt0}

\begin{figure}[h]
$$\resizebox{120mm}{!}{\includegraphics{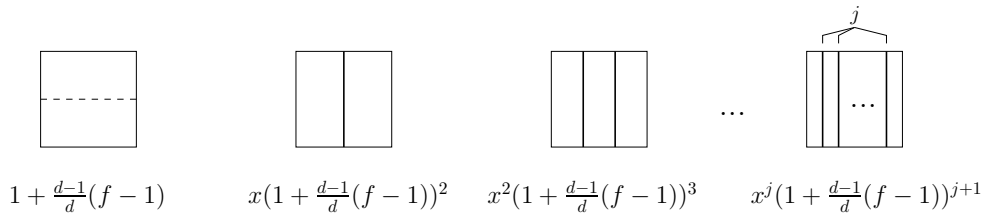}}$$
\caption{Illustration to the proof of Theorem~\ref{the:alt0}}
\label{fig:alt0}
\end{figure}

\noindent \textit{Proof.}
Let $j$ be the number of principal cuts of a partition $P$, orthogonal to $x_1$-axis, $0\leq j<m$.
They split $B$ into $j+1$ parts all of which may have all possible $m$-alternating partitions with
a principal cut in another direction.
Therefore partitions with exactly $j$ such cuts contribute $x^{j}( 1+\frac{d-1}{d}(f-1))^{j+1}$ to the generating function, see Fig.~\ref{fig:alt0}. Therefore
\[f=\sum_{j=0}^{m-1} x^{j}\left( 1+\frac{d-1}{d}(f-1)\right)^{j+1} =
\left(1 + \frac{d-1}{d}(f-1)\right)
\left(
\frac
{1 - x^m \left( 1+\frac{d-1}{d}(f-1) \right)^m}
{1 - x \left( 1+\frac{d-1}{d}(f-1) \right)}
\right) .
\]

Denote $g=1 + \frac{d-1}{d}(f-1)$ -- this is the generation function enumerating $m$-alternating partitions with a fixed direction of the principal cut. We have

\[\frac{d}{d-1}(g-1) + 1 = g \cdot \left( \frac{1 - x^m g^m}{1-xg}\right) \]
which gives

\[ g = 1+(d-1)xg^2 \cdot \left( \frac{1-x^{m-1}g^{m-1}}{1-xg} \right) .\]
$\hfill\square$

\vspace{2mm}

\begin{figure}[ht]
$$\resizebox{80mm}{!}{\includegraphics{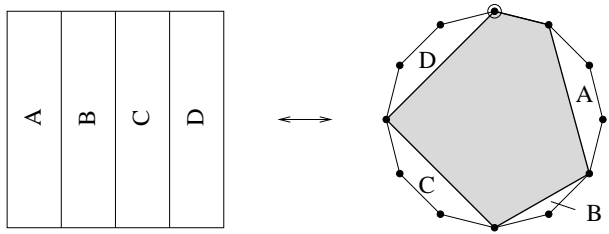}}$$
\caption{Illustration to Remark 3}
\label{fig:poly}
\end{figure}

We consider some special cases.
\begin{enumerate}
\item If $m$ is not bounded, we obtain all the guillotine partitions. Indeed, in this case we have
\[f = \left(1 + \frac{d-1}{d}(f-1)\right)
\frac
{1}
{1 - x \left( 1+\frac{d-1}{d}(f-1) \right)},\]
which gives after simplifications
\[f = 1+xf+(d-1)xf^2.\]

\item Let $m=2$. In this case we have
$g = 1+(d-1)xg^2$.
Since $h = 1 + xh^2$ is the generating function for Catalan numbers, we have
\[f = 1 + \sum_{k\geq 1}d(d-1)^{k-1}C_kx^k,\]
or
\[a_d(n) = d(d-1)^{n-1}C_n\]
for $n \geq 1$.
This can be also easily proved by induction.

\item Let $d=2$. We have
\[ g = 1+xg^2 \cdot \left( \frac{1-x^{m-1}g^{m-1}}{1-xg} \right) .\]

This $g$ is also known to be the generating function enumerating the number of dissections of a convex polygon with $(n+2)$ vertices by non-crossing diagonals into polygons with at most $m+1$ vertices. Indeed, it is easy to construct a bijection between $m$-alternating guillotine partitions with $n$ cuts, principal cut being in a fixed direction, and such dissections. A subbox with $j$ principal cuts would correspond to a $j+2$-gon in the dissection, see Fig.~\ref{fig:poly} (the rotated letters in the left side mean that non-trivial ``inner'' partitions have horizontal principal cuts).
\end{enumerate}

\vspace{2mm}

We obtain an explicit formula for the number of $m$-alternating partitions of a $d$-box by $n$ cuts. If we denote $h = \frac{f-1}{d}$ (where $f$ is the generating function from Theorem~\ref{the:alt0}), we obtain

\[h = \sum_{i=2}^{m} x^{i-1} \left(1+(d-1)h \right)^i= \frac { x(1+(d-1)h)^2 - x^m (1+(d-1)h)^{m+1} } { 1 - x(1+(d-1)h) }. \]

We use Lagrange's inversion formula (see \cite[Section 5.4]{S}) to obtain:

\[h = \sum_{p \geq 1} \frac{1}{p} [h^{p-1}] \left(\frac {x^p (1+(d-1)h)^{2p} (1-x^{m-1}(1+(d-1)h)^{m-1})}
 {(1 - x(1+(d-1)h))^p}\right )^p = \]

\[= \sum_{p \geq 0} \frac{x^{p+1}}{p+1} [h^{p}]
\left(
\frac
{(1+(d-1)h)^{2p+2} \sum_{i=0}^{p+1}(-1)^i \binom{p+1}{i} x^{(m-1)i} (1+(d-1)h)^{(m-1)i}}
{(1 - x(1+(d-1)h))^{p+1}}
\right) =
\]

\[
= \sum_{p\geq 0} \sum_{i=0}^{p+1} \sum_{j\geq 0}
\frac{(-1)^i}{p+1} \binom{p+1}{i} \binom{p+j}{j} \binom{2p+2+(m-1)i+j}{p}
(d-1)^p
x^{p+1+(m-1)i+j}.
\]

Taking $n=p+1+(m-1)i+j$ we finally obtain that
the number of $m$-alternating partitions of a $d$-box by $n$ cuts ($n \geq 1$) is

\[a_d(n) = d \sum_{p\geq 0} \sum_{i=0}^{p+1}
\frac{(-1)^i}{p+1} \binom{p+1}{i} \binom{n-1-(m-1)i}{p} \binom{n+p+1}{p}
(d-1)^p.\]

\subsection{Guillotine partitions avoiding \protect\winaa}
\label{sec:win}

Suppose $p$ is a principal cut of a subbox of $B$,
and two
parallel cuts $q$ and $r$ belong to different subspaces bounded
by $p$ and meet $p$.
The authors of~\cite{ack} do not distinguish the cases where
$q$ is higher than $r$ and conversely (\ \winaa vs. \winbb), and they mention
the enumeration when such partitions are considered different as an open problem.
We were able to enumerate partitions which \emph{avoid} this situation.
Thus, we consider partitions of $B$ with the following
restriction: No subbox of $B$ has a principal cut $p$ and two
parallel cuts $q$ and $r$ that belong to different subspaces bounded
by $p$ and meet $p$. We shall refer to these partitions as to \emph{partitions
avoiding \winaa}.

\begin{figure}[ht]
$$\resizebox{120mm}{!}{\includegraphics{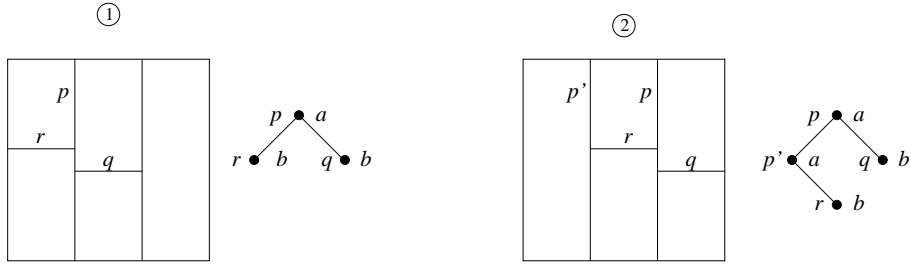}}$$
\caption{Two types of subtrees which correspond to \protect\winaa}
\label{fig:win1}
\end{figure}

Assume without loss of generality that $q$ belongs to
the upper subspace bounded by $p$, and $r$ belongs to the lower
subspace bounded by $p$.
There are two cases.
If $p$ is the lowest principal cut in its subbox, then $r$ is its left child, $q$ is its right child. Otherwise, if $p'$ is the next principal cut below $p$, then $p'$ is the left child of $p$, $r$ is the right child of $p'$, $q$ is the right child of $p$. See Figure~\ref{fig:win1}.

\vspace{2mm}

Let $f = f_{d, a}$ be the generating
function counting non-empty trees with the root having a fixed color $a$.
Then we have (See Fig.~\ref{fig:win4}):
\begin{figure}[ht]
$$\resizebox{135mm}{!}{\includegraphics{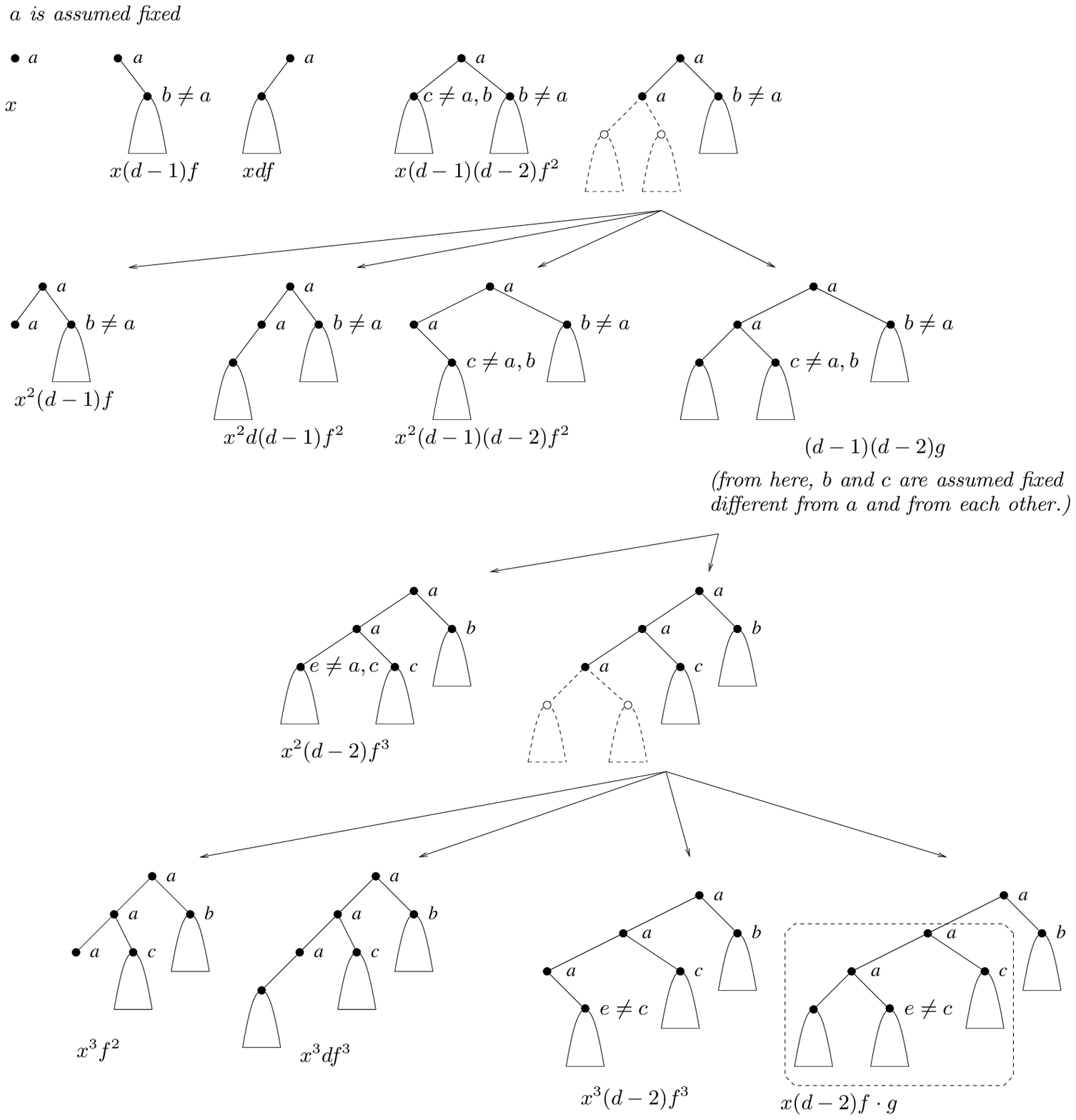}}$$
\caption{Trees corresponding to partitions which avoid \protect\winaa} \label{fig:win4}
\end{figure}

$$\begin{array}{l} f = x + x(d-1)f + xdf + x(d-1)(d-2)f^2 + x^2(d-1)f\\
\qquad\qquad\qquad\qquad\qquad+ x^2d(d-1)f^2 + x^2(d-1)(d-2)f^2 +
(d-1)(d-2)g,\end{array}$$ where $g$ satisfies

\[g = x^2(d-2)f^3 + x^3f^2 + x^3df^3 + x^3(d-2)f^3 + x(d-2)f\cdot g.\]

This gives

\[g = \frac{x^2(d-2)f^3 + x^3f^2 + x^3df^3 + x^3(d-2)f^3}{1-x(d-2)f}.\]

Thus

$$\hspace{-58mm}f = x + x(d-1)f + xdf + x(d-1)(d-2)f^2 + x^2(d-1)f +$$
$$\hspace{5mm}+ x^2d(d-1)f^2 +x^2(d-1)(d-2)f^2+ (d-1)(d-2)\dfrac{x^2(d-2)f^3 +
x^3f^2 + x^3df^3 + x^3(d-2)f^3}{1-x(d-2)f},$$

which implies

\[f =x+x(x-1+2d)f+dx(x+d-2)f^2\]

Hence,

\[f=\frac{1+x-x^2-2dx-\sqrt{(1+x-x^2-2dx)^2-4dx^2(x+d-2)}}{2dx(x+d-2)}=\]
\[=\frac{x}{1+x-x^2-2dx}C\left(\frac{dx^2(x+d-2)}{(1+x-x^2-2dx)^2}\right),\]
where $C(t)=\frac{1-\sqrt{1-4t}}{2t}$ is the generating function for
the Catalan numbers $c_n=\frac{1}{n+1}\binom{2n}{n}$. In order to
get an explicit formula for the $x^n$ coefficient of $f$ we recall that
$C(t)=\sum_{n\geq0}\frac{1}{n+1}\binom{2n}{n}t^n$, and obtain

\[ \hspace{-67mm} f=\sum\limits_{k\geq0}c_k\dfrac{d^kx^{2k+1}(x+d-2)^k}{(1+x-x^2-2dx)^{2k+1}}= \]
 \[ \hspace{-11mm} =\sum\limits_{k\geq0}\sum\limits_{j\geq0}\sum\limits_{i=0}^kc_k d^k(d-2)^{k-i} \binom{k}{i}\binom{2k+j}{j} x^{2k+i+j+1}(x+2d-1)^j = \]
\[=\sum\limits_{k\geq0}\sum\limits_{j\geq0}\sum\limits_{i=0}^k\sum_{\ell=0}^j
c_k d^k(d-2)^{k-i}(2d-1)^{j-\ell}
\binom{k}{i}\binom{2k+j}{j}\binom{j}{\ell} x^{2k+i+j+1+\ell},
\]
which implies that the $x^n$ coefficient of $f$ is

$$\sum\limits_{k=0}^{(n-1)/2}\sum\limits_{j=0}^{n-1-2k}\sum\limits_{i=0}^k c_k d^k(d-2)^{k-i}(2d-1)^{2j+2k+i+1-n}
\binom{k}{i}\binom{2k+j}{j}\binom{j}{n-2k-i-j-1}.$$ Hence, we can
state the following result.

\begin{wind}
\label{the:wind}
The number of Guillotine partitions of a
$d$-dimensional box by $n$, $n\geq1$, cuts avoiding \winaa is
$$d\sum\limits_{k=0}^{(n-1)/2}\sum\limits_{j=0}^{n-1-2k}\sum\limits_{i=0}^k c_k d^k(d-2)^{k-i}(2d-1)^{2j+2k+i+1-n}
\binom{k}{i}\binom{2k+j}{j}\binom{j}{n-2k-i-j-1}.$$
\end{wind}

\vspace{2mm}

Consider the planar case ($d=2$) of the following variation: now the forbidden situation is that
$p$ is a principal cut of a subbox of $B$,
and two parallel cuts $q$ and $r$ belong to different subspaces bounded
by $p$ (but they do not necessarily meet $p$). In other words, we consider \wincc-avoiding partitions.
Let $f$ be the generating function counting such partitions, $h$ the generating function counting all such non-trivial partitions with principal cut in a fixed direction. Then we have

\[h = \sum_{k\geq 1} (1+(k+1)h)x^k \]
which gives
\[h = \frac{x}{1-x}+h\cdot\left( \frac{1}{(1-x)^2}-1 \right),\]
\[h = \frac{x(1-x)}{2(1-x)^2-1},\]
and finally
\[f = \frac{1-2x}{1-4x+2x^2}.\]

\subsection{Restricted guillotine partitions and permutations}
\label{sec:restr_par_per}

In view of the planar case of Theorem~\ref{the:permd}, we can interpret the results on restricted guillotine partitions, with $d=2$, in terms of permutations. Figure~\ref{fig:perm_part} shows how patterns in guillotine partitions correspond to patterns in separable permutations: (1) presents the pattern forbidden in boundary partitions, (2) in \winaa-avoiding partitions, (3) in \wincc-avoiding partitions. We obtain the following result:

\begin{figure}[ht]
$$\resizebox{150mm}{!}{\includegraphics{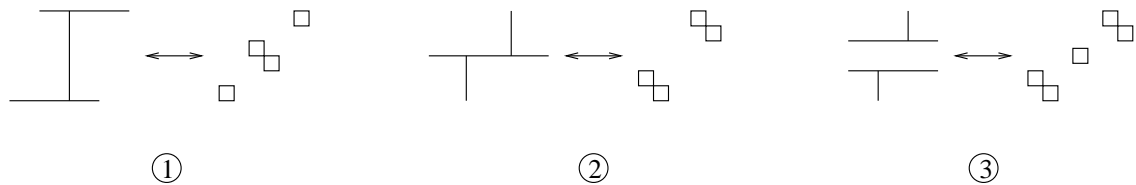}}$$
\caption{Patterns in guillotine partitions and the corresponding patterns in permutations} \label{fig:perm_part}
\end{figure}

\begin{enumerate}

  \item Permutations avoiding $2413$, $3142$, $1324$, $4231$ are in the correspondence with planar guillotine boundary partitions. The generating function of their enumerating sequence is
      $f=1+\frac{2x(1-3x+3x^2)^2}{(1-x)(1-2x)^4}$, and the first ten terms are $1$, $2$, $6$, $20$, $64$, $194$, $562$, $1570$, $4258$, $11266$
(the first column ($d=2$) in Table~\ref{tab:boundary}).

\item Permutations avoiding $2413$, $3142$, $21\bar354$, $45\bar312$ correspond to \winaa-avoiding partitions. The generating function of their enumerating sequence satisfies $f =x+x(x+3)f+2x^2f^2$,
    its first ten numbers are $1$, $2$, $6$, $20$, $70$, $254$, $948$, $3618$, $14058$, $55432$. This sequence is \cite[A078482]{S} --
     it is also mentioned by Atkinson and Stitt it~\cite{AS} as the counting sequence of number of permutations of $[n]$ avoiding $2413, 3142, 1342, 1423$. \emph{(Can a bijection between these two families of restricted permutations be constructed?)}

  \item Permutations avoiding $2413$, $3142$, $2143$, $3412$ correspond to \wincc-avoiding guillotine partitions. The generating function of their enumerating sequence is $f = \frac{1-2x}{1-4x+2x^2}$ and the first ten terms $1$, $2$, $6$, $20$, $68$, $232$, $792$, $2704$, $9232$, $31520$. This sequence is
      \cite[A006012]{S}.

\end{enumerate}
\bibliographystyle{amsplain}

\end{document}